\documentclass[12pt]{article}
\usepackage{amsmath,amssymb}

\newcommand{\al}{\hspace{1em}}
\newcommand{\ad}{\al \em}

\newtheorem{thm}{\bf \al Theorem}
\newtheorem{prop}{\bf \al Proposition}
\newtheorem{lem}{\bf \al Lemma}

\newtheorem{defi}{\bf \al Definition}

\newcommand{\Pf}{\al Proof. \ad}
\newcommand{\qd}{\hfill $\mathbf \square$}

\begin{document}
\renewcommand{\textfraction}{0}


\title{Results on zeta functions for codes
\footnote{
Presented at the Fifth Conference on Algebraic Geometry, Number Theory,
Coding Theory and Cryptography, University of Tokyo, January 17-19, 2003}
}
\author{Iwan Duursma\footnote{Supported by NSF Grant DMS-0099761.
Address: Department of Mathematics, University of Illinois 
at Urbana-Champaign, 1409 W. Green Street, Urbana IL 61801, USA.
E-mail: duursma@math.uiuc.edu.}
}
\date{February 05, 2003}
\maketitle

\thispagestyle{empty}


\begin{abstract}
We give a new and short proof of the Mallows-Sloane upper bound 
for self-dual codes. We formulate a version of Greene's theorem
for normalized weight enumerators. We relate normalized 
rank-generating polynomials to two-variable zeta functions.
And we show that a self-dual code has the Clifford property, 
but that the same property does not hold in general for formally 
self-dual codes.
\end{abstract}

\section{Introduction} \label{sec:int}

In \cite{Duu1} we introduced, for an arbitrary linear code, its zeta function,
as a different way to describe the weight distribution of the code. 
The definition is motivated by properties of algebraic curves and of codes 
constructed with those curves. After analyzing the definition more carefully
for its coding theoretic meaning, we formulated in \cite{Duu2} an equivalent 
definition in terms of puncturing and shortening operations. Both definitions
are recalled in this paper together with some basic properties of zeta 
functions for codes. \\
We introduce a polynomial $g(w)$ of small degree that 
interpolates the normalized differences
\[
\left( \frac{A_w}{{n \choose w}} - (q-1)\frac{A_{w-1}}{{n \choose w-1}}
\right) (-1)^{w-d},
\]
for $w=1,2,\ldots,n$ (Lemma \ref{lem:gw}). 
The fact that the polynomial is both of small degree
and has many zeros when the minimum distance $d$ is large leads us to
an alternative proof for the Mallows-Sloane upper bounds
(as a special case of the bounds in Theorem \ref{thm:DC}). 
The polynomial $g(w)$ determines the zeta polynomial of 
a linear code and vice versa (Proposition \ref{prop:gwPT}). \\
Pellikaan defined a two-variable zeta function for curves \cite{Pel}.  
For codes, we can consider a similar two-variable zeta function. 
In the approach that we take here, we first formulate a version
of Greene's Theorem for normalized rank-generating polynomials
(\ref{eq:AnWn}). Then we define the two-variable zeta function
in terms of the normalized rank-generating polynomial 
(\ref{eq:ZTWn}). 
And we show that this is compatible with Definition \ref{def:PT2} for 
the one-variable zeta function. \\
Clifford's theorem on the dimension of special divisors has an analogue
for codes. We use an argument from \cite{Oxl} to show that the corresponding
result holds for self-dual codes, 
but in general not for formally self-dual codes. 
 
\hyphenation{mo-di-fied}
\hyphenation{enu-me-ra-tor}

\section{Weight enumerators and zeta functions} \label{sec:we}

Let $C$ be a linear code of length $n$ and minimum distance $d$
over the finite field of $q$ elements.
Let $A_i$ be the number of words of weight $i$ in $C$. 
The weight enumerator of the code $C$ is defined as
\[
A(x,y) = x^n + \sum_{i=d}^n A_i x^{n-i}y^i
\]
\begin{defi}[\cite{Duu1}] \label{def:PT1}
\ad For a given weight enumerator $A(x,y)$, of a $q$-ary linear code
of length $n$ and minimum distance $d$, define $P(T)$ as the unique
polynomial of degree at most $n-d$ such that     
\[
[T^{n-d}]~ \frac{P(T)}{(1-T)(1-qT)} \; (y+(x-y)T)^n = \frac{A(x,y)-x^n}{q-1} 
\]
\end{defi}
Let $a_w = A_w / {n \choose w}$, for $w=0,1,\ldots,n$. Define the
normalized weight enumerator as
\[
a(t) = \frac{1}{q-1} (a_d + a_{d+1}t + \cdots + a_n t^{n-d})
\]
\begin{defi}[\cite{Duu2}] \label{def:PT2}
\ad For a given normalized weight enumerator $a(t)$, of a $q$-ary linear code
of length $n$ and minimum distance $d$, define $P(T)$ as the unique
polynomial of degree at most $n-d$ such that  
\[
\frac{P(T)}{(1-T)(1-qT)}(1-T)^{d+1} ~ \equiv 
   ~ a(\frac{T}{1-T}) \pmod{T^{n-d+1}}
\]
\end{defi}
As a brief motivation for Definition \ref{def:PT1}, consider the special case 
$P(T)=1$, and recall that 
\[
\frac{1}{(1-T)(1-qT)}
\]
is a generating function for the number of monic polynomials of degree at most
a given degree $a$, say. To interpret the modified generating function in the 
definition, we use
\[
\frac{(y(1-T)+xT)}{1-T} = y + xT + xT^2 + \cdots
\]
It is then clear that the coefficient at $x^{n-i} y^i T^{a}$ 
gives the number of those monic polynomials of degree at most $a$ that have 
precisely $n-i$ zeros in a given subset $\{x_1, \ldots, x_n\} \subset F_q$.
Thus the weight enumerator $A(x,y)$ that corresponds to $P(T)=1$ is realized
by the linear code 
\[
C = \{ (f(x_1), \ldots, f(x-n) : f \in F_q[x]_{\leq a} \}
\]
The code $C$ has $k=a+1, d=n-a$ and meets the Singleton bound $d \leq n-(k-1)$. 
Definition \ref{def:PT2} is motivated by the following property of the
normalized weight enumerator.

\begin{thm}[\cite{Duu2}] \label{thm:atinv} 
The expression 
\[
a(t) (1+t)^d \pmod{t^{n-d+1}}
\]
is invariant under puncturing or shortening. 
\end{thm}

To have well-defined puncturing (projection) 
and shortening (restriction) operations on a weight enumerator, 
independent of the choice of a coordinate, we average over all coordinates, 
so that the effect on the weight enumerator $A(x,y)$ is given by
\newcommand{\pr}{\partial}
\[
\frac{1}{n}( \frac{\pr}{\pr_x} + \frac{\pr}{\pr_y} )  ~~\text{(puncturing)}
\quad \quad
\frac{1}{n}( \frac{\pr}{\pr_x} ) ~~\text{(shortening)}
\]
The following properties are derived in \cite{Duu1}.
For nondegenerate linear codes, with both $d \geq 2$ and $d^\perp \geq 2$, 
\[
\deg P(T) = n+2-d-d^\perp \quad \text{ and } \quad P(1) = 1 
\]
Duality, as contained in the MacWilliams identities
\[
A_{C^\perp}(x,y) = \frac{1}{|C|}A_C(x+(q-1)y,x-y),
\]
becomes
\[
P^\perp(T) = P(1/qT) \: q^g T^{g+g^\perp}, 
\]
where $g=n+1-k-d$ and $g^\perp=n+1-k^\perp-d^\perp$. \\

The zeros of the zeta polynomial play a r\^ole in the following upper bound
for the minimum distance. 
Writing $P(T)=a_d/(q-1)(1+aT+\cdots)$, Definition \ref{def:PT2} yields
\[
a_d (a - d + q) = a_{d+1}
\]
or 
\[
d+1 = q+1+a - a_{d+1}/a_d \leq q+1+a
\]
Thus, estimates for the reciprocal zeros of $P(T)$, and in particular for 
their sum $-a$, yield upper bounds for the minimum distance of
a linear code. \\

The following theorem describes the zeros of $P(T)$ for an
interesting infinite family of weight enumerators. For a self-dual code of
type (IV), that is to say defined over $F_4$ with only words of even weight,  
$3d \leq n + 6$ (\cite{MacSlo}). 
When the bound is met the weight
enumerator of the code is uniquely determined.
 
\begin{thm}[\cite{Duu4}]
Let $A(x,y)$ be the unique weight enumerator of type (IV) with $d=m+3$ and
$n=3m+3$, for $m$ odd. Let $P(T)$ be the associated zeta polynomial and
let $Q(T) = P(T)(1+2T)$. Then
\[
Q(\frac{T^2}{2}) = \lambda_m C_m^{m+1} (\frac{T^{-1}+T}{2}) T^m
\]
where $C_m^{m+1}$ is an ultraspherical polynomial of degree $m$ with
$m$ real zeros on $[-1,1]$, and $\lambda_m$ a constant depending on $m$. 
In particular $Q(e^{i2\theta}/2)=0$ if and only if $C_m^{m+1}(\cos \theta)=0$.
\end{thm}

Details for this section and related results can be found in 
\cite{Duu1}, \cite{Duu2}, \cite{Duu3}, \cite{Duu4}.
 
\section{The Mallows-Sloane bounds}
 
As in the previous section, let $a_w = A_w / {n \choose w}$.

\begin{lem} \label{lem:gw}
For a linear code of length $n$, minimum distance $d$ and dual
minimum distance $d^\perp$,
there exists a unique polynomial $g(w)$ of degree $n-d^\perp$ 
such that
\[
g(w) = \left( a_w - (q-1)a_{w-1} \right) (-1)^{w-d}, 
\quad \text{ for }w=1,2,\ldots,n. 
\]
\Pf We choose $g(w)$ of degree at most $n-d^\perp$ such that it 
interpolates the right hand side correctly for $w=1,2,\ldots,n-d^\perp+1$. 
It remains to show
(1) $g(w)$ interpolates correctly in $n-d^\perp+2, \ldots, n$, and
(2) $\deg g(w) = n-d^\perp$.
Let $S \subset \{1,2,\ldots,n\}$ be a subset of size $s$ and consider the
subcode of $C$ of words with support on $S$. Averaging over all $S$ of
size $s$ gives for the average size of such a subcode
\[
\sum_{w=0}^n a_w {s \choose w}
\]
On the other hand for $s > n-d^\perp$, the size of each such subcode equals
\[
q^{k-(n-s)}
\]
Thus, for $s > n-d^\perp$,
\begin{align*}
&\sum_{w=0}^{s+1} a_w {s+1 \choose w} = 
    q \sum_{w=0}^{s} a_w {s \choose w}, \\
&\sum_{w=1}^{s+1} a_w {s \choose w-1} = 
    (q-1) \sum_{w=0}^{s} a_w {s \choose w}, \\
&\sum_{w=0}^{s} {s \choose w} (a_{w+1} - (q-1) a_{w}) = 0.
\end{align*}
With elementary calculus, this says that the value for
\[
(-1)^w (a_{w} - (q-1) a_{w-1})
\]
at $w=s+1$ is the polynomial extrapolation of the values at $w=1,\ldots,s$.
This proves claim (1). For $s=n-d^\perp$, the average size of a subcode
exceeds $q^{k-d^\perp}$, and the extrapolation relation cannot be used
to obtain the value at $w=n-d^\perp+1$ from the values at 
$w=1,\ldots,n-d^\perp$. This clearly implies claim (2). \qd
\end{lem}

Since $g(w)$ has zeros at $2,3,\ldots,d-1$, we obtain $d-2 \leq n-d^\perp$.
And when equality holds,
\[
g(w) = (q-1) {w-2 \choose d-2} (-1)^{w-d}
\]
For a general weight enumerator, let $P(T) = p_0 + p_1 T + \cdots + P_r T^r$
be the zeta polynomial, with $r=n+2-d-d^\perp$. Then

\begin{prop} \label{prop:gwPT}
\[
g(w) = (q-1) \left( p_0 {w-2 \choose d-2} - p_1 {w-2 \choose d-1} + \cdots +
                    (-1)^r p_r {w-2 \choose d+r-2} \right)
\]
\end{prop}

\begin{thm}[\cite{Duu4}] \label{thm:DC}
Let the code $C$ have all weights divisible by $c$. Then
\[
d+cd^\perp \leq n+c(c+1)
\]
If moreover the code is binary, even, and contains the allone word, then 
\[
2d+cd^\perp \leq n+c(c+2)
\]
\Pf We give a proof based on Lemma \ref{lem:gw}.
From $g(w)$ we can obtain a polynomial $h(w)$ of same degree such that
\[ 
h(w) = \left( a_{w} - (q-1)^c a_{w-c} \right) (-1)^{w-d}, 
\quad \text{ for }w=c,c+1,\ldots,n. 
\]
It has at least $(c-1)/c \cdot (n-c) + 1/c \cdot (d-2c)$ zeros, and
\[
(c-1)n-(c-1)c+d-2c \leq cn-cd^\perp ~\Leftrightarrow~ 
d+cd^\perp \leq n+c(c+1)
\]
For the second claim, the degree of $h(w)$ drops to at most $n-d^\perp-1$. 
It has at least $(c-1)/c \cdot (n-c) + 2/c \cdot (d-2c)$ zeros, and
\[
(c-1)n-(c-1)c+2d-4c \leq cn-cd^\perp-c ~\Leftrightarrow~ 
2d+cd^\perp \leq n+c(c+2)
\]
\qd
\end{thm}
When applied to self-dual codes, with $d=d^\perp$, we recover the
Mallows-Sloane upper bounds (\cite{MacSlo}).
\begin{align*} \label{eq:masl}
\begin{array}{llll}
\text{\rm Type I}   &(q=2,c=2): & &d \leq 2 \lfloor n/8 \rfloor + 2. \\
\text{\rm Type II}  &(q=2,c=4): & &d \leq 4 \lfloor n/24 \rfloor + 4. \\
\text{\rm Type III} &(q=3,c=3): & &d \leq 3 \lfloor n/12 \rfloor + 3. \\
\text{\rm Type IV}  &(q=4,c=2): & &d \leq 2 \lfloor n/6 \rfloor +2.  
\end{array}
\end{align*}

\section{Two-variable zeta functions}

Pellikaan defined, for an algebraic curve over a finite field,
the two-variable zeta function as the convergent power series 
\[
Z(T,u) = \sum_{[D]} \frac{u^{l(D)}-1}{u-1} T^{\deg D}
\]
The summation is over divisor classes $[D]$. 
For a finite field of size $q$ and for $u=q$, it agrees with the Hasse-Weil 
zeta function: $Z(T,q) = Z(T)$. 
Some familiar properties of the Hasse-Weil zeta function generalize to the
Pellikaan zeta function \cite{Pel}. Thus $Z(T,u)$ is a 
rational function in the variables $T$ and $u$, with functional equation
\[
Z(T,u) = Z(\frac{1}{uT},u) u^{g-1} T^{2g-2}
\]
The vanderGeer-Schoof two-variable zeta function gives a generalization
to number fields. In the version for curves it is defined as 
\[
\zeta^{GS}(s,t) =\sum_{[D]} q^{s h^0 + t h^1}
\]
where $h^0 = \dim L(D)$ and $h^1 = \dim \Omega(D) = \dim L(K-D)$.
We use it in the form  
\[
Z^{GS}(x,y) = \sum_{[D]} x^{h_0} y^{h_1}
\]
so that $\zeta^{GS}(s,t) = Z^{GS}(q^{s},q^{t})$.
Deninger gives the relation between $Z(T,u)$ and $\zeta^{GS}(s,t)$ 
(Proposition 2.1 \cite{Den}). For $Z^{GS}$ it becomes,
\begin{equation} \label{eq:Den}
Z(T,u) (u-1) T^{1-g} = Z^{GS}(uT,T^{-1})
\end{equation}

Two-variable rank-generating polynomials for matroids go back to Whitney 
and to important papers in graph theory by Tutte. The columns in the 
generating matrix of a code form a set $G$. For each subset $A$ of columns,
let
\begin{align*} 
r(A) &= {\rm{rank}}(A)  &\text{  (rank)}\\
\rho(A) &= |A|    &\text{  (degree)}\\
n(A) &= \rho(A)-r(A)  &\text{  (nullity)} 
\end{align*}
The rank-generating polynomial (or Whitney polynomial, or corank-nullity 
polynomial, e.g. \cite{BryOxl}) is defined as
\[
W_G(x,y) = \sum_{A \in G} x^{r(G)-r(A)} y^{|A|-r(A)}
\]
For a code with column set $G$, the weight enumerator is given by
Greene's Theorem \cite{Gre}, which can be written in the form
\[
\frac{A(x,y)}{(x-y)^k y^{n-k}} = W_G( \frac{qy}{x-y}, \frac{x-y}{y})
\]
The rank-generating polynomial of a code depends only on the generators of 
the code. To compute the weight enumerator of a code after taking 
coefficients in an extension field, only $q$ needs to be replaced. 
We give a version of Greene's theorem for the normalized rank-generating
polynomial. Let
\begin{align}
W_n(x,y) &= \sum_{i=0}^n \frac{1}{{n \choose i}} 
      \sum_{A \in G, |A|=i} x^{r(G)-r(A)} y^{|A|-r(A)}. \label{eq:Wn} \\
A_n(x,y) &= \sum_{i=0}^n \frac{1}{{n \choose i}} A_i x^{n-i} y^i \label{eq:An}
\end{align}
Then 
\begin{align} 
A_n(s,t)(s+t)^{n+1} = 
  ~~&W_n(\frac{qt}{s+t},\frac{s+t}{t})(s+t)^k t^{n-k} s^{n+1} 
         \nonumber \\
  +&\tilde W_n(\frac{qs}{s+t},\frac{s+t}{s})(s+t)^k s^{n-k} t^{n+1}
         \label{eq:AnWn}
\end{align}
The relation is written as a polynomial identity but the polynomial 
$\tilde W_n$ has a priori no particular meaning, so the relation could 
as well be used with $s=1$ as a congruence relation modulo $t^{n+1}$,
\begin{equation} \label{eq:AnWnt}
A_n(1,t)(1+t)^{n+1} \equiv 
  W_n(\frac{qt}{1+t},\frac{1+t}{t})(1+t)^k t^{n-k} \pmod{t^{n+1}}
\end{equation}
An exception is for binary self-complementary codes that have 
$A_n(s,t)=A_n(t,s)$ and $\tilde W_n = W_n$. 
We want to show that there is a natural definition of a two-variable
zeta function for codes that is compatible with our earlier definitions
for the one-variable case. \\

To relate the two-variable zeta function of Pellikaan and the 
rank-generating polynomial, let, for a special divisor $E$ on the curve,
\begin{align*}
r(E) &= l(K)-l(K-E) \\
\rho(E) &= \deg(E) \\ 
n(E) &= \deg(E)-(l(K)-l(K-E)) = l(E)-1
\end{align*}
These definitions do not make the set of special divisors into a
representable matroid (unless we allow a somewhat wider definition)
but they seem perfectly natural and give a satisfactory correspondence.
The canonical divisor $K$ has rank $l(K)-1$ and the rank-generating
polynomial for special divisors becomes
\[
W(x,y) = \sum_{[E]} x^{l(K-E)-1} y^{l(E)-1}
\]
which is similar to the vanderGeer-Schoof two-variable zeta function. 

To define a two-variable zeta function for codes, we use two
properties of the two-variable zeta function for curves:
(1) the number of divisor classes of given degree is constant
and equal to $h$.
(2) the zeta function consists of a finite contribution and an 
infinite tail that only depends on $h$.
The first property holds with $h=1$ for the normalized rank-generating
polynomial $W_n(x,y)$.
For the second property we add an infinite tail to $W_n$.
\begin{equation} \label{eq:Wn+}
W_n^+(x,y) = W_n(x,y) + \frac{x^{k+1}}{1-x} + \frac{y^{n-k+1}}{1-y}.
\end{equation}
For a normalized rank-generating function $W_n^+$, define a
two-variable zeta function, in analogy with (\ref{eq:Den}), via 
\begin{equation} \label{eq:ZTWn}
Z(T,u)(u-1) T^{1-g} = W_n^+(uT,T^{-1})
\end{equation}
We show that this definition is compatable with the one-variable zeta 
function in Definition \ref{def:PT2}, such that $Z(T,q)=Z(T)$. 
We modify (\ref{eq:AnWnt}) to include contributions
of the infinite tail that was added to $W_n$ in (\ref{eq:Wn+}).
Let $x=qt/(1+t)$.
\[
(\frac{qt}{1+t})^{k+1} \frac{1+t}{1-(q-1)t} (1+t)^k t^{n-k}
\equiv 0 \pmod{t^{n+1}}
\]
Let $y=(1+t)/t$.
\[
(\frac{1+t}{t})^{n-k+1} \frac{t}{-1} (1+t)^k t^{n-k}
\equiv -(1+t)^n \pmod{t^{n+1}}
\]
Combined with (\ref{eq:AnWnt}) and (\ref{eq:Wn+}) this gives
\[
(A_n(1,t)-1)(1+t)^{n+1} 
  \equiv W_n^+(\frac{qt}{1+t},\frac{1+t}{t}) 
                   (1+t)^k t^{n-k} \pmod{t^{n+1}}.
\]
Or, using $A_n(1,t)-1 = (q-1)a(t)t^d$ and (\ref{eq:ZTWn}),
\[
a(t)t^d (1+t)^{n+1} \equiv Z(\frac{t}{1+t},q) (\frac{t}{1+t})^{1-g} 
                                 (1+t)^k t^{n-k} \pmod{t^{n+1}}.
\]
Finally, with $d=n+1-k-g$ this reduces to
\[
a(t) (1+t)^{d+1} \equiv Z(\frac{t}{1+t},q) \pmod{t^{n+1-d}}
\]
which agrees with Definition \ref{def:PT2} after the substitution
$t=T/(1-T)$. Compare also with Theorem \ref{thm:atinv}. \\

As an example, an MDS code of length $n$ and dimension $k$         
has
\[
W_n(x,y) = x^k + \cdots + x + 1 + y + \cdots + y^{n-k},
\]
\[
W_n^+(x,y) = \frac{1-xy}{(1-x)(1-y)},
\]
and
\[
Z(T,u) = \frac{-T^{-1}}{(1-uT)(1-T^{-1})} = \frac{1}{(1-T)(1-uT)}   
\]
The passage from $W_n$ to $W_n^+$ to define the two-variable zeta 
function of a code is in line with Theorem \ref{thm:atinv}. The effect of 
puncturing on the polynomial $W_n$ is $W_n \mapsto W_n - y^{n-k}$,
and the effect of shortening is $W_n \mapsto W_n - x^{k}$. Thus by
adding an infinite tail, $W_n^+$ has become invariant under
puncturing or shortening.

\section{A Clifford type theorem for self-dual codes}

We give an interpretation of Clifford's theorem for self-dual codes. 
In \cite{TsfVla}, \cite{Vla}, Clifford's theorem is used to give
estimates for the weight distributions of geometric Goppa codes.
Clifford's theorem says
\[ 
l(E)-1 + l(K-E)-1 \leq l(K)-1
\]
which corresponds to an inequality 
\[
n(A) + r(G) - r(A) \leq r(G)  \Leftrightarrow  2 r(A) \geq |A|
\]
for representable matroids. 

\begin{prop}
The inequality $2 r(A) \geq |A|$ holds for any code that contains
its dual and for any choice of columns $A$. In the special case of
a self-dual code $C$, equality holds if and only if $C=C_1 \oplus C_2$ 
for self-dual codes $C_1$ and $C_2$ that are supported on $A$ and the 
complement of $A$, respectively. \\

\Pf Corank and nullity are dual notions and 
$k^\perp-r^\perp(\bar A)=|A|-r(A)$. 
The subcode of the dual code with support on $A$ therefore has dimension 
$|A|-r(A)$. The subcode is self-orthogonal and thus $2|A|-2r(A)\leq |A|$, 
with equality if and only if it is self-dual. 
Using duality twice, we have 
$k^\perp+|\bar A|-2r^\perp(\bar A) = k+|A|-2r(A)$. And for a self-dual
code, $2r(A) = |A|$ if and only if $2r(\bar A) = |\bar A|$ if and only
if both $A$ and $\bar A$ support selfdual codes, in which case clearly
$C= C_1 \oplus C_2$ as required. \qd
\end{prop}

A short argument to prove the Clifford inequality for self-dual codes
is provided by \cite[Theorem 3.9]{Oxl}.
The $n=2k$ columns in a self-dual code divide in at least one 
way into two independent subsets of size $k$ each. Let the subset
$A$ have $a_1$ columns in the first subset and $a_2$ columns in the
second subset. Then $2 r(A) \geq 2 \max \{a_1,a_2\} \geq a_1+a_2 = |A|$. \\

The inequality $2 r(A) \geq |A|$ for self-dual codes, does in general not 
hold for formally self-dual codes. It is easy to find a formally self-dual 
code for which the inequality fails. 
We may take $(1 0)$ with dual code $(0 1)$. 

\bibliographystyle{plain}

\end{document}